\documentclass{article}
\usepackage{amsfonts, amssymb, latexsym, amscd}
\usepackage{amsmath, amsthm}
\newtheorem{thm}{Theorem}[section]

\def\G{{\Gamma}}
\def\setunion{{ \cup }}

\begin{document}

\title{\bf Cospectral pairs of regular graphs with different connectivity}

\author{Willem H. Haemers\thanks{e-mail haemers@uvt.nl}
\\
{\it\small Department of Econometrics and Operations Research,}
\\
{\it\small Tilburg University, The Netherlands}
}
\date{}

\maketitle

\begin{abstract}
\noindent
For vertex- and edge-connectivity we construct infinitely many pairs of regular graphs
with the same spectrum, but with different connectivity.
\\[5pt]
{Keywords:}~graph spectrum, vertex-connectivity, edge-connectivity, spectral characterization.
AMS subject classification:~05C50.
\end{abstract}

\begin{center}
{\em Didicated to the memory of Slobodan Simi\'{c}}
\end{center}

\section{Introduction}

Spectral graph theory deals with the relation between the structure of a graph and the eigenvalues (spectrum) of an associated matrix,
such as the adjacency matrix $A$ and the Laplacian matrix $L$.
Important types of relations are the spectral characterization.
These are conditions in terms of the spectrum of $A$ or $L$, which are necessary and sufficient for certain graph properties.
Two famous examples are: (i) a graph is bipartite if and only if the spectrum of $A$ is invariant under multiplication by $-1$, and
(ii) the number of connected components of a graph is equal to the multiplicity of the eigenvalue $0$ of $L$.
Properties that are characterized by the spectrum for $A$ as well as for $L$
are the number of vertices, the number of edges, and regularity.
If a graph is regular, the spectrum of $A$ follows from the spectrum of $L$, and vice versa.
This implies that for both $A$ and $L$ the properties of being regular and bipartite,
and being regular and connected are characterized by the spectrum.

If a property is not characterized by the spectrum, then there exist a pair of cospectral graphs
where one has the property and the other one not.
For many graph properties and several types of associated matrices, such pairs are not hard to find.
However, if we restrict to regular graphs it becomes harder and more interesting,
because a pair of regular cospectral graphs where one has a given property and the other one not
is a counter examples for a spectral characterization with respect to $A$, $L$, and several other types of matrices.
Such a pair of regular cospectral graphs
has been found for a number of properties, for example
for being distance-regular~\cite{H}, having a given diameter~\cite{HS}, and admitting a perfect matching~\cite{BCH}.

The {\em vertex-connectivity}~$\kappa(\G)$ of a graph $\G$ is the minimum number of vertices one has to delete from $\G$ such that the graph
becomes disconnected.
The {\em edge-connectivity}~$\kappa'(\G)$ is the minimum number of edges one has to delete from $\G$ to make the graph disconnected.
One easily has that $\kappa(\G)\leq\kappa'(\G)\leq\delta(\G)$ where $\delta(\G)$ is the minimal degree of $\G$.
Clearly $\kappa(\G)=0$ as well as $\kappa'(\G)=0$ just means that $\G$ is disconnected,
therefore these two properties are characterized by the spectrum when $\G$ is regular.
Fiedler~\cite{F} showed that the second smallest eigenvalue of the Laplacian matrix $L$ (called the {\em algebraic connectivity})
is a lower bound for the vertex- (and edge-) connectivity.
For a regular graph $\G$ there exist stronger spectral bounds for $\kappa(\G)$ (see \cite{A}) and $\kappa'(\G)$ (see~\cite{C}).
Here we show that for the vertex- and for the edge-connectivity in a connected regular graph there is in general no spectral characterization.
For $k \geq 2$ we present a pair of regular cospectral graphs $\G$ and $\G'$ of degree $2k$ and order $6k$,
where $\kappa(\G)=2k$ and $\kappa(\G')=k+1$.
The edge-connectivity turned out to be much harder.
Nevertheless, for every even $k\geq 4$ we found a pair of regular cospectral graphs $\G$ and $\G'$ of degree
$3k-5$, where $\kappa'(\G)=3k-5$ and $\kappa'(\G')=3k-6$.

The main tool is the following result of Godsil and McKay~\cite{GM}.
\begin{thm}\label{gm}
Let $\G$ be a graph, and let $X_1,\ldots, X_m,Y$ be a partition of the vertex set of $\G$
into $m+1$ classes, such that the following holds:\\
(i) For $1\leq i,j\leq m$, each vertex $x\in X_i$ has the same numbers of neighbors in $X_j$.
\\
(ii) For $1\leq i\leq m$ each vertex $y\in Y$ is adjacent to $0$, $\frac{1}{2}|X_i|$, or all vertices of $X_i$.
\\
Make a graph $\G'$ as follows:
For $0\leq i\leq m$ and each vertex $y\in Y$ with $\frac{1}{2}|X_i|$ neighbors in $X_i$,
delete the $\frac{1}{2}|X_i|$ edges between $y$ and $X_i$, and insert $\frac{1}{2}|X_i|$
edges between $y$ and the other vertices of $X_i$.
Then the graph $\G'$ thus obtained is cospectral with $\G$ with respect to the adjacency matrix.
\end{thm}
The corresponding operation is called {\em Godsil-McKay switching}.
In many applications $m=1$.
Then $X_1$ is called a {\em GM-switching set}, and condition (i) just means that the subgraph of $\G$ induced by $X_1$ is regular.

We assume familiarity with basic results from linear algebra and graph spectra; see for example \cite{BH}, or~\cite{CRS}.
As usual, $J_{m,n}$ (or just $J$) denotes the $m\times n$ all-ones matrix, $O_{m,n}$ (or just $O$) is the $m\times n$ all-zeros matrix,
and $I_n$ (or just $I$) denotes the identity matrix of order~$n$.

\section{Vertex-connectivity}

Suppose $k\geq 2$.
We define a $k$-regular graph $G$ on the integers modulo $3k-1$ as follows:
for $i=0,\ldots,2k-1$ vertex $i$ is adjacent to $\{k+i,k+i+1,\ldots,2k+i-1\}$ (mod~$3k-1$).
Then $G$ has no triangles and vertex-connectivity $k$ (indeed, between any pair of vertices there exists
$k$ vertex-disjoint paths).
Next we partition the vertex set $V$ of $G$ into four classes $V_0,\ldots,V_3$ as follows:
\[
V_0=\{0\},\ V_1=\{k,k+1,\ldots,2k-1\},\ V_2=\{1,2,\ldots,k-1\},\ V_3=\{2k,2k+1,\ldots,3k-2\}.
\]
So $V_1$ consists of the neighbors of vertex~$0$.
Note that $G$ contains a matching of size $k-1$ that matches vertices of $V_2$ with $V_3$.
Let $B$ be the corresponding partitioned adjacency matrix of $G$.
Then
\[
\ B=\left[
\begin{array}{cccc}
0 & J_{1,k} & O_{1,k-1} & O_{1,k-1} \\
J_{k,1} & O_{k,k} & B_{1,2} & B_{1,3} \\
O_{k-1,1} & B_{1,2}^\top & O_{k-1,k-1} & B_{2,3} \\
O_{k-1,1} & B_{1,3}^\top & B_{2,3}^\top & O_{k-1,k-1}
\end{array}
\right].
\]
Next we define
\[
\ N=\left[
\begin{array}{c}
I_{k+1} \\
J_{k-1,k+1}
\end{array}
\right],
\ K =
\left[
\begin{array}{c}
O_{k,k-1} \\
J_{k,k-1}
\end{array}
\right],
\ M =
\left[\, J\!-\!N\ \ K\ \ J\!-\!K\, \right].
\]
Finally, we define $\G$ to be the graph with adjacency matrix
\[
A=\left[
\begin{array}{ccc}
O_{2k,2k} & N                & M \\
N^\top    & \ J\!-\!I_{k+1}\ & O_{k+1,3k-1} \\
M^\top    & O_{3k-1,k+1}     & B
\end {array}
\right].
\]
Let $\{X,U,V\}$ be the corresponding partition of the vertex set of $\G$.
It follows that $\G$ is regular of degree $2k$, and that $X$ is a GM-switching set of $\G$.
In terms of $A$, Godsil-McKay switching replaces $N$ by $J-N$ and $M$ by $J-M$.
The graph $\G'$ thus obtained is cospectral with $\G$, and becomes disconnected if we delete the first $k+1$ vertices.
So the vertex-connectivity of $\G'$ is at most $k+1$.

To verify that the vertex-connectivity of $\G$ equals $2k$ we have to find $2k$
vertex-disjoint paths between any two distinct nonadjacent vertices $x$ and $y$ of $\G$ (see for example \cite{S}, Theorem~15.1).
This is a straightforward (and time consuming) activity, for which we have to distinguish several cases,
depending on the partition classes to which $x$ and $y$ belong.

Suppose $x,y\in V$. Then, because $G$ has vertex-connectivity $k$, there exist $k$ vertex-disjoint paths in $G$ between $x$ and $y$.
Let $X_x$ and $X_y$ be the sets of neighbors in $X$ of $x$ an $y$, respectively.
Then there exist $\ell = |X_x\cap X_y|$ vertex-disjoint paths between $x$ and $y$ of length $2$,
and $k-\ell$ vertex-disjoint paths between $X_x\setminus X_y$ and $X_y\setminus X_x$ via $U$ of length $3$.

Suppose $x,y\in X$.
If $x$ and $y$ are both adjacent to all vertices of $U$, then $x$ and $y$ are also adjacent to all vertices of $V_3$,
so there are $2k$ vertex-disjoint paths of length $2$ between $x$ and $y$.
If $x$ and $y$ are both adjacent to all vertices of $V_2$, then $x$ and $y$ have $2k-2$ common neighbors in $V$, and there is a
path between $x$ and $y$ of length $3$ via two vertices of $U$, and a path of length $4$ via a vertex of $X$.
If $x$ is adjacent to all vertices of $U\setunion V_3$, and $y$ is adjacent to all vertices of $V_2$,
then there exist $k-1$ vertex-disjoint paths of length $4$ using the matching between $V_2$ and $V_3$,
$x$ and $y$ have one common neighbor $z\in U$ and $k$ vertex-disjoint paths of length $5$ via $V_0\setunion V_1$, $X$ and $U\setminus\{z\}$.
Next suppose $x$ is the unique vertex in $X$ adjacent to all vertices of $V_3$ and just one vertex of $U$.
If $y$ is adjacent to all vertices of $U\setunion V_3$, then $x$ and $y$ have one common neighbor in $U$ and $k-1$ common neighbors in $V_3$,
furthermore there are $k$ vertex-disjoint paths of length $4$ via $V_1$, $X$ and $U$.
If $y$ is adjacent to all vertices of $V_2$, then $x$ and $y$ have $k-1$ vertex-disjoint paths of length $2$ via $V_1$,
$k-1$ disjoint paths of length $3$ via $V_3$ and $V_2$, one path of length $3$ via $V_0$ and $V_1$, and one path of length $3$ via an edge of $U$.

The remaining cases: $(x\in X,~y\in U)$, $(x\in X,~y\in V)$ and $(x\in U,~ y\in V)$ are left as an exercise.
In a similar way one can verify that $\G'$ has vertex connectivity $k+1$.
So we can conclude:

\begin{thm}
For every $k\geq 2$ there exists a pair of $2k$-regular cospectral graphs,
where one has vertex-connectivity $2k$ and the other one has vertex-connectivity $k+1$.
\end{thm}

The smallest cases ($k=2,3,4$) have been double checked by computer using the package newGRAPH~\cite{SBCS}.

\section{Edge-connectivity}\label{ec}

The construction of cospectral pairs of regular graphs with different edge-connectivity turned out to be much harder then
for the case of vertex-connectivity.
Again Godsil-McKay switching is the main tool, but now we apply Theorem~\ref{gm} with $m=2$.
(All attempts to find an example with $m=1$ failed; the difficulty is caused by the regularity requirement.)
For every even integer $k\geq 6$ we define a graph $\G$ for which the vertex set is partitioned into three classes, $X_1$, $X_2$ and $Y$,
and assume the corresponding partitioned adjacency matrix $A$ has the following structure.
\[
A=\left[
\begin{array}{ccc}
A_1 & L & M_1\\
L^\top & A_2 & M_2\\
M_1^\top & M_2^\top & B
\end{array}
\right],
\]
such that $A_1$, $A_2$ and $B$ are the adjacency matrices of the subgraphs $G_1$, $G_2$ and $H$ induced by $X_1$, $X_2$,
and $Y$, respectively.
Let $G$ be the graph induced by $X=X_1\cup X_2$, defined by
\[
L = \left[\begin{array}{cc} J_{k-1,k-1} & O \\ O & J_{k+1,k+1}-{C} \end{array}\right], \mbox{ and }
A_1=A_2=L-I_{2k},
\]
where $C$ is the adjacency matrix of the $(k+1)$-cycle.
Moreover, $H$ is a disconnected graph with adjacency matrix
\[
B = \left[\begin{array}{cccc}
B_1 & O & J_{2k-3,k-1} & O \\
O & B_2 & O & J_{2k-5,k+1} \\
J_{k-1,2k-3} & O & J-I_{k-1} & O\\
O & J_{k+1,2k-5} & O & J-I_{k+1}
\end{array}\right],
\]
where $B_1$ is the adjacency matrix of a $(k-4)$-regular graph $H_1$ of order $2k-3$,
and $B_2$ is the adjacency matrix of a $(k-6)$-regular graph $H_2$ of order $2k-5$
(here we use that $k\geq 6$ and even).
Finally we define
\[
\left[\begin{array}{c}
M_1 \\ M_2
\end{array}\right] =
\left[\begin{array}{cccccc}
O & J_{k,k-2} & O & O & O_{k,2k} \\
O & O & J_{k,k-2} & O & O_{k,2k}  \\
J_{k,k-2} & O & O & O & O_{k,2k}  \\
O & O & O & J_{k,k-2} & O_{k,2k}
\end{array}\right].
\]
A vertex from $H$, which is not a vertex of $H_1$ or $H_2$ is adjacent to no vertex of $X$.
A vertex from $H_1$ or $H_2$ is adjacent to exactly half or none of the vertices of $X_1$ and $X_2$.
It is easily checked that $G_1$, $G_2$ and $G$ are regular,
thus we can conclude that the given partition of $\G$ satisfies the conditions for GM-switching
given in Theorem~\ref{gm}.
To obtain the adjacency matrix $A'$ of the switched graph $\G'$ we have to replace $M_1$ and $M_2$ by
\[
\left[\begin{array}{c}
M'_1 \\ M'_2
\end{array}\right] \mbox{ by }
\left[\begin{array}{cccccc}
O & O         & J_{k,k-2} & O & O_{k,2k} \\
O & J_{k,k-2} & O         & O & O_{k,2k} \\
O         & O & O & J_{k,k-2} & O_{k,2k} \\
J_{k,k-2} & O & O & O         & O_{k,2k}
\end{array}\right].
\]
We know that $\G$ and $\G'$ are cospectral, and we easily have that $\G$ and $\G'$ are ($3k-5$)-regular of order $10k-8$.

To work out the edge-connectivity of $\G$ and $\G' $, we first observe that the graph $G_1$ (resp. $G_2$) has two
connected components $G_{1,1}$ and $G_{1,2}$ (resp. $G_{2,1}$ and $G_{2,2}$),
and we define $x_1$ (resp. $x_2$) to be the first vertex (ordered as in $A$) of the larger component $G_{1,2}$ (resp. $G_{2,2}$).
Let $y$ be the last vertex of $H_1$.
For $i,j=1,2$ let $X_{i,j}$ be the vertex set of $G_{i,j}$.
Then the adjacencies between $X$ and $Y$ are as follows.
\\
(i) the first $k-2$ vertices of $B_1$ are adjacent to $X_{2,1}\setunion\{x_2\}$ in $\G$, and to $X_{2,2}\setminus\{x_2\}$ in $\G'$;
\\
(ii) the second $k-2$ vertices of $B_1$ are adjacent to  $X_{1,1}\setunion \{x_1\}$ in $\G$, and to $X_{1,2}\setminus\{x_1\}$ in $\G'$;
\\
(iii) the vertex $y$ of $B_1$ is adjacent to $X_{1,2}\setminus\{x_1\}$ in $\G$, and to $X_{1,1}\setunion\{x_1\}$ in $\G'$;
\\
(iv) the first $k-1$ vertices of $B_2$ are adjacent to $X_{1,2}\setminus\{x_1\}$ in $\G$, and to $X_{1,1}\setunion\{x_1\}$ in $\G'$;
\\
(v) the other $k-2$ vertices of $B_2$ are adjacent to $X_{2,2}\setminus\{x_2\}$ in $\G$, and to $X_{2,1}\setunion\{x_2\}$ in $\G'$.
From (i) to (v) we see that $\G$ and $\G'$ become disconnected if we delete the vertices $x_1$, $x_2$ and $y$
($\kappa(\G)=\kappa(\G')=3$).
We claim that there is no disconnecting set of edges in $\G$, which is smaller than the degree $3k-5$.
The only candidates for such a set are subsets of the edges between $X$ and $Y$.
The smallest disconnecting set of edges between $X$ and $Y$ has $3k-4$ edges and consists of the $2k-4$ edges
between $\{x_1,x_2\}$ and $H_1$, together with the $k$ edges between $y$ and $X_{1,2}$.
Therefore $\G$ has edge-connectivity $3k-5$.
However, after switching, we find a disconnecting edge set in $\G'$ consisting of the $2k-5$ edges between $\{x_1,x_2\}$ and $H_2$,
and the $k-1$ edges between $y$ and $X_{1,1}$ (indeed, we don't need the edge $\{y,x_1\}$).
Therefore the edge-connectivity of $\G'$ equals $3k-6$.
Thus we have:
\begin{thm}
For every $k\geq 6$ there exists a pair of $(3k-5)$-regular cospectral graphs, where one has edge-connectivity $3k-5$
and the other one has edge-connectivity $3k-6$.
\end{thm}
The smallest pair has degree 13, order 52 and edge-connectivities 13 and 12.
There is some variation possible in the above construction, which can lead to examples with other
degrees, orders and edge-connectivities.
For example, we can obtain a pair of 7-regular graphs of order 36 with edge connectivities 7 and 6 by the above construction with $k=4$
when we replace the component of $H$ containing $H_2$ by a graph with adjacency matrix
\[
\left[
\begin{array}{cccc}
O_{3,3} & I_3   & I_3   & I_3   \\
I_3 & J-I_3 & J-I_3 & J-I_3 \\
I_3 & J-I_3 & J-I_3 & J-I_3 \\
I_3 & J-I_3 & J-I_3 & J-I_3
\end{array}
\right].
\]
But we found no pair of cospectral regular graphs where one has edge-connectivity smaller than 6.
The cases $k=4$ and $k=6$ have been double checked by computer using the package newGRAPH~\cite{SBCS}.

\section{Final remarks}

If $\G$ and $\G'$ are regular cospectral graphs, then also the line graphs $L(\G)$ and $L(\G')$
are regular and cospectral (see for example~\cite{BH}, Section~1.4.5).
Moreover, for every graph $\G$, $\kappa(L(\G)) = \kappa'(\G)$.
So, the line graphs of the graphs described in Section~\ref{ec} give another infinite family of
cospectral pairs of regular graphs with different vertex-connectivity.

We can conclude that in general the property of being regular with a given vertex- or edge-connectivity
is not characterized by the spectrum.
However, for vertex-connectivity at most 2, and edge-connectivity at most 5 this may still be the case.
Especially interesting is the question if being regular with vertex- or edge-connectivity 1 is characterized
by the spectrum.
\\[10pt]
{\bf Acknowledgment}
I thank the referees for pushing me to be more explicit in the construction of Section~2;
it made the construction better understandable.

\end{document}